\documentclass[11pt]{article}
\usepackage[]{geometry}
\usepackage[mathcal]{euscript}

\usepackage{bm}                     
\usepackage{amsfonts}
\usepackage{amssymb}
\usepackage{amsmath}
\usepackage{amsthm}
\usepackage{graphicx}               
\usepackage{pstricks,pst-plot}
\usepackage{wrapfig}
\usepackage{colortbl}               
\usepackage{booktabs}

\usepackage{amssymb,amsmath}

\usepackage[latin1]{inputenc}
\usepackage{graphicx,color}
\newtheorem{definition}{Definition}
\newtheorem{remark}{Remark}
\newtheorem{theorem}{Theorem}
\newtheorem{lemma}{Lemma}
\newtheorem{proposition}{Proposition}

\begin{document}

\begin{center}
\large \bf  On statistical properties of sets fulfilling rolling-type conditions \normalsize
\end{center}
\normalsize

\

\begin{center}
  Antonio Cuevas$^{a}$, Ricardo Fraiman$^{b}$, Beatriz Pateiro-L\'opez$^{c}$ \\
  $^{a}$ Universidad Aut\'onoma de Madrid, Spain\\
  $^{b}$ Universidad de San Andr\'es, Argentina and   Universidad de la Rep\'ublica, Uruguay\\
  $^{c}$ Universidad de Santiago de Compostela, Spain\\
\end{center}

\begin{abstract}
Motivated by set estimation problems, we consider three closely related shape conditions for compact sets: positive reach, $r$-convexity and rolling condition.
First, the relations between these shape conditions are analyzed.
Second, we obtain for the estimation of sets fulfilling a rolling condition a result of ``full consistency'' (i.e., consistency with respect to the Hausdorff metric for the target set and for its boundary).
Third, the class of uniformly bounded compact sets whose reach is not smaller than a given constant $r$ is shown to be a $P$-uniformity class (in Billingsley and Tops{\o}e's (1967) sense) and, in particular, a Glivenko-Cantelli class.
Fourth, under broad conditions, the $r$-convex hull of the sample is proved to be a fully consistent estimator of an $r$-convex support in the two-dimensional case. Moreover, its boundary length is shown to converge (a.s.) to that of the underlying support.
Fifth, the above results are applied to get new consistency statements for level set estimators based on the excess mass methodology (Polonik, 1995).
\end{abstract}



\section{Introduction}\label{sec:intro}

Three geometric closely related properties, with fairly intuitive interpretations,  are analyzed in this paper. They are called \sl positive reach\/\rm, (see, e.g., Federer \cite{federer:59}, Rataj \cite{rataj:05}, Ambrosio et al. \cite{ambrosio:08}), \sl $r$-convexity,\/ \rm (Perkal \cite{perkal:56}, Mani-Levitska \cite{mani:93}, Walther \cite{walther:97}) and  \sl rolling condition\/ \rm (Walther \cite{walther:97, walther:99}). Our interest in these ``rolling-type properties'' is motivated in the framework of set estimation, that is, the problem of reconstructing a set $S\subset{\mathbb R}^d$ (typically a density support or a density level set) from a random sample of points; see e.g. Cuevas and Fraiman \cite{cuevas:09} for a recent survey. The classical theory of this subject, dating back to the 1960's, is largely concerned with the assumption that $S$ is convex; see e.g. D{\"u}mbgen and Walther \cite{dumbgen:96} or Reitzner \cite{reitzner:09} for a survey. The rolling-type properties have been employed as  shape restrictions on $S$ alternative to (and much broader than) the convexity assumption.

\

The rolling-type conditions are useful in statistics and stochastic
geometry at least in two ways. First, they can be sometimes incorporated to the estimator: for example, if the (compact) support $S$ of a random variable $X$ is assumed to be $r$-convex one could estimate $S$, from a random sample of $X$ by taking the $r$-convex hull of the sample points, much in the same way as the convex hull has been used to estimate a convex support; see Rodr\'{i}guez-Casal \cite{rodriguez:07}. Second, they can be used as regularity assumptions in order to get faster rates of convergence for the estimators; see, e.g., Cuevas and Rodr\'{i}guez-Casal \cite{cuevas:04}, Pateiro-L\'{o}pez and Rodr\'{i}guez-Casal \cite{pateiro:08}.
Of course other, perhaps more standard, regularity conditions have also been used in set estimation. They rely on usual smoothness assumptions on the boundary or the underlying density, defined in terms of derivatives. See Biau et al. \cite{biau:09, biau:08} and Mason and Polonik \cite{mason:09} for some recent interesting examples.
Some  deep results on the connection between differentiability assumptions and rolling-type conditions can be found  in Federer \cite{federer:59}, Walther \cite{walther:99} and Ambrosio et al. \cite{ambrosio:08}.

The most popular among the rolling-type properties is by far the positive reach condition. It was introduced by Federer \cite{federer:59} in a celebrated paper which could be considered as a  landmark in geometric measure theory. Among other relevant results, Federer \cite{federer:59} proved that (Theorem 5.6), for a compact set $S$ with reach $r>0$, the volume (Lebesgue measure) of the $\epsilon$-parallel set $B(S,\epsilon)$ can be expressed
as a polynomial in $\epsilon$, of degree $d$, for $0\leq \epsilon<r$. This is a partial generalization of the classical Steiner formula that shows this property for convex sets for all $\epsilon>0$. So the positive reach property can be seen as a natural generalization of convexity in a much deeper sense than that suggested by the definition. For some recent interesting contributions on this property see Ambrosio et al. \cite{ambrosio:08} and Colesanti and Manselli \cite{colesanti:10}. An application in set estimation, more specifically in the problem of estimating the boundary measure, can be found in Cuevas et al. \cite{cuevas:07}.

The $r$-convexity property provides a different but closely related generalization of convexity (see Section \ref{rolling} for precise definitions). Whereas this property has also a sound intuitive motivation it is much less popular. An earlier reference is Perkal \cite{perkal:56} but, to our knowledge, the first statistical application is due to Walther \cite{walther:97} who uses this condition in the setting of level set estimation. A study of the $r$-convex hull as an estimator of an $r$-convex support can be found in Rodr\'{i}guez-Casal \cite{rodriguez:07}.

\

Let us know establish some notation and basic definitions. We
are concerned here with subsets of ${\mathbb R}^d$ although some concepts and results can be stated, with little additional effort, in the broader setup of metric spaces. The Euclidean norm in ${\mathbb R}^d$ will be denoted by $\Vert\cdot\Vert$.
Given a set
$A\subset {\mathbb R}^d$, we denote by $A^c$, $\textnormal{int}(A)$ and $\partial A$
the complement, interior and boundary of $A$, respectively. We denote by
$B(x,r)$  the closed ball with
centre $x$ and radius $r$. By convenience, the open ball $\mbox{int}(B(x,r))$ will be denoted by $\mathring{B}(x,r)$.

In the problem of estimating a compact set and/or its boundary we need to use some suitable distances in order to assess the quality of the estimation and establish asymptotic results (concerning consistency, convergence rates and asymptotic distribution). The most usual distances in this setting are the \sl Hausdorff distance\/ \rm $d_H$   and the \sl distance in measure\/ \rm $d_\nu$. The definitions are as follows.

Let $\mathcal{M}$ be the
class of closed, bounded, nonempty subsets of ${\mathbb R}^d$. For
$A,C\subset\mathcal{M}$, the Hausdorff distance between $A$ and
$C$ is defined by
\begin{equation*}\label{haus2}
d_H(A,C)=\inf\left\{\epsilon>0:\ \
C\subset B(A,\epsilon) {\textnormal{ and }}
A\subset B(C,\epsilon)\right\},
\end{equation*}
where $B(A,\epsilon)$ denotes the closed
$\epsilon$-parallel set of $A$,
$B(A,\epsilon)=\{x\in {\mathbb R}^d:\ \delta_A(x)\leq\epsilon\}$, with
$\delta_A(x)=\inf\{\Vert x-y\Vert:\, y\in A\}$. $({\mathcal M},d_H)$ is a complete locally compact metric space.

Let $\nu$ be a Borel measure on ${\mathbb R}^d$, with $\nu(C)<\infty$ for any compact $C$. Let $A,C$ be Borel sets with finite $\nu$-measure. The distance in measure between $A$ and $C$ is defined by
\begin{equation*}\label{mdist}
d_\nu(A,C)=\nu(A\Delta C),
\end{equation*}
where $\Delta$ denotes the symmetric difference between $A$ and $C$, that is, $A\Delta C=(A\setminus C)\cup(C\setminus A)$. Often $\nu$ is either a probability measure or   the Lebesgue measure on $\mathbb{R}^d$, which we will denote by $\mu$. Note that the distance function $d_\nu$ is actually a pseudometric but it becomes a true metric if we identify two sets differing in a $\nu$-null set. This amounts to work in the quotient space associated with the corresponding equivalence relation.

Besides the estimation of the set $S$ and the boundary $\partial S$, some functionals of $S$ are sometimes of interest as estimation targets. This is the case of the  $(d-1)$-dimensional boundary measure, $L(S)$ (that is  the perimeter of $S$ for $d=2$ and the surface area for $d=3$). There are several, not always equivalent, definitions for $L(S)$ (see Mattila \cite{mattila:95}) but we will mainly use the \sl outer Minkowski content\rm,
\begin{equation}\label{os}
L(S)=\lim_{\epsilon\to 0}\frac{\mu(B(S,\epsilon)\setminus S)}{\epsilon}.
\end{equation}
Under regularity conditions, the limit in (\ref{os}) coincides with the usual \sl Minkowski content\/ \rm given by
\begin{equation}\label{mc}
L_0(S)=\lim_{\epsilon\to 0}\frac{\mu(B(\partial S,\epsilon))}{2\epsilon}.
\end{equation}
We refer to Ambrosio et al. \cite{ambrosio:08} for general conditions ensuring the existence of the outer Minkowski content and its relation to other measurements of the boundary of $S$.

\

The contributions in this paper can be now summarized as follows. In Section \ref{rolling} we clarify the relations between the rolling-type conditions pointing out that if the reach of a set is $r$, then it is $r$-convex and, in turn, $r$-convexity entails the $r$-rolling property. See Propositions \ref{reachrconvex} and \ref{rconvexroll}. It is also shown that the converse implications are not true in general. 
In Section \ref{sub:seq} we prove that, under very general conditions on $\nu$ and $S$, $d_H(S_n,S)\to 0$ plus $d_H(\partial S_n,\partial S)\to 0$ (we call this simultaneous convergence \sl full convergence\rm)  implies $d_\nu(S_n,S)\to 0$ (Theorem \ref{dHdmu}). We also show (Theorem \ref{partialSn}) that if the sets $S_n$ fulfill the rolling condition, the $d_H$-convergence $d_H(S_n,S)\to 0$ implies full convergence. This result (which has some independent interest) is used below in the paper.
In Section \ref{sec:bill} we show that, under broad conditions, the class of uniformly bounded sets with reach $\geq r$ is a  $P$-uniformity class (in the sense of Billingsley and Tops{\o}e \cite{billingsley:67}) and therefore also a Glivenko-Cantelli class. This is interesting from the point of view of empirical processes, see e.g., Devroye et al. \cite{devroye:96} or van der Vaart \cite{vaart:98}, and will be also used in Section \ref{sec:appl} of this paper.
In subsection \ref{sub:reach} we show that, if a support $S$ in ${\mathbb R}^2$ is assumed to be $r$-convex then, the natural estimator (which, as mentioned above, is the $r$-convex hull of the sample) is consistent in all the usual senses. In particular, it provides also a plug-in consistent estimator of the boundary length whose practical performance is checked through some numerical comparisons in the appendix. This is an interesting contribution to the theory of nonparametric boundary estimation which so far relies mostly on the use of two samples (one inside and the other outside the set $S$); see Cuevas et al. \cite{cuevas:07}, Pateiro-L\'{o}pez and Rodr\'{i}guez-Casal \cite{pateiro:08} and Jim\'{e}nez and Yukich \cite{jimenez:10}.
The results mentioned in the above points are used in subsection \ref{EM} for the problem of estimating density level sets of type $\{f\geq \lambda\}$ using the excess mass approach (see Polonik \cite{polonik:95}). In particular, we obtain new consistency properties for the estimation of $\{f\geq \lambda\}$ as well as uniform consistency results (with respect to $\lambda$) for the same problem.

\section{Rolling-type assumptions}\label{rolling}
Convexity is a natural geometrical restriction that arises in many
fields. While the study of this topic dates back to antiquity, most important contributions and applications date from the 19th and 20th centuries. We refer to the Handbook edited by Gruber and Wills \cite{gruber:93} for a complete survey of convex geometry and its relations to other areas of mathematics.  However, convexity
leaves out usual features of sets such as holes or notches and may
be considered an unrealistic assumption in some situations. As a natural consequence, there
are meaningful extensions of the notion of convexity; see e.g. Mani-Levitska \cite{mani:93}.
In this Section, we shall consider three different shape restrictions that generalize that of convexity. We group them under the name of rolling-type conditions. Their formal definitions are as follows.



\begin{definition}
Given $r>0$,
a set $S\subset\mathbb{R}^d$ is said to fulfill the (outside)
$r$-rolling condition if for all $x\in \partial S$ there is a
closed ball with radius $r$, $B_x$, such that $x\in B_x$ and
$\mbox{int}(B_x) \cap S=\emptyset$.
\end{definition}
The radius $r$ of the ball acts a smoothing parameter. We refer to
Walther \cite{walther:99} for a deep study of this condition.


\

Following the notation
by Federer \cite{federer:59}, let $\textnormal{Unp}(S)$ be the set of
points $x\in\mathbb{R}^d$ having a unique projection on $S$, denoted by $\xi_S(x)$. That is, for $x\in \textnormal{Unp}(S)$, $\xi_S(x)$ is the unique point which minimizes $\delta_S(x)$.

\begin{definition}
For $x\in S$ let
$\textnormal{reach}(S,x) = \sup\{r>0:\ \mathring{B}(x,r)\subset\textnormal{Unp}(S)\}.$
The reach of $S$ is then defined by
\begin{equation*}
\textnormal{reach}(S) = \inf\{\textnormal{reach}(S,x):\ {x\in S}\}
\end{equation*}
and $S$ is said to be of positive reach if $\textnormal{reach}(S)
> 0$. 
\end{definition}
Besides all convex closed sets (which have infinite reach) and
regular submanifolds of class 2 of $\mathbb{R}^d$, the class of
sets of positive reach  also contains nonconvex sets or sets whose
boundary is not a smooth manifold. Sets with positive reach were introduced by Federer \cite{federer:59} who also obtained their main properties.
In particular, these sets obey a Steiner formula in the following sense.

\begin{theorem}\label{fed5.6} (Federer (1959, Th. 5.6))
Let $S\subset\mathbb{R}^d$ a compact set with $\textnormal{reach}(S)>0$. Let $K$ be a Borel subset of $\mathbb{R}^d$. Then there exist unique Radon measures $\Phi_0(S,\cdot)$, $\Phi_1(S,\cdot),\ldots, \Phi_d(S,\cdot)$ over $\mathbb{R}^d$ such that, for $0\leq \epsilon<\textnormal{reach}(S)$,
\begin{equation}\label{steiner_d}
\mu(B(S,\epsilon)\cap\{x:\ \xi_S(x)\in K\})=\sum_{i=0}^d{\epsilon^{d-i}b_{d-i}\Phi_i(S, K)}
\end{equation}
where $b_0=1$ and, for $j\geq 1$, $b_j$ is the $j$-dimensional measure of a unit ball in $\mathbb{R}^j$.
\end{theorem}

\begin{remark}\label{remmc}
Theorem \ref{fed5.6} plays an important role in some problems of stochastic geometry regarding the calculation of boundary measurements.
The measures $\Phi_j:=\Phi_1(S,\cdot)$ are the so-called \it curvature measures \/ \rm associated with $S$, for $j=0,1,\ldots,d$.
In particular, if $\textnormal{reach}(S)>0$, it is straightforward from Theorem \ref{fed5.6} with $K=S$, that the outer Minkowski content $L(S)$  is finite, and corresponds to the first order coefficient in the expansion (\ref{steiner_d}). Furthermore, by Corollary 3 and inequality (27) in Ambrosio et al. \cite{ambrosio:08}, combined with Federer's Theorem, we have that if $\textnormal{reach}(S)>0$, then  $\partial S$ has finite Minkowski content $L_0(S)$, which in particular yields $\mu(\partial S)=0$. For this reason we do not impose in what follows $\mu(\partial S)=0$ whenever we assume positive reach for $S$.
\end{remark}


\begin{definition}
A set $S\subset \mathbb{R}^d$ is said to be
$r$-convex for some $r>0$ if $S=C_r(S)$, where
\begin{equation}\label{rhull}
  C_r(S)=\bigcap_{\mathring{B}(y,r)\cap
    S=\emptyset}\mathring{B}(y,r)^c.
\end{equation}
\end{definition}
That is, $S$ is
$r$-convex if for all $x\in S^c$, there exists
$y\in\mathring{B}(x,r)$ such that $\mathring{B}(y,r)\cap
S=\emptyset$. We refer to Perkal \cite{perkal:56} for elementary properties of $r$-convex sets and connections between convexity and $r$-convexity. Walther \cite{walther:99} and Rodr\'{i}guez-Casal \cite{rodriguez:07} also deal with this shape restriction in the context of set estimation. A natural question is whether certain characterizations of convex sets are still meaningful in the context of generalized notions of convexity. For instance, given a set $S$ and $r>0$ it is possible to find the minimal $r$-convex set containing $S$, the so-called \it $r$-convex hull \/\rm of $S$. In fact, it follows from the properties of $r$-convex sets that the $r$-convex hull of $S$ coincides with the set $C_r(S)$ given in (\ref{rhull}), see Perkal \cite{perkal:56}. Note that the definition of $C_r(S)$ resembles that of the convex hull (with balls of radius $r$ instead of halfspaces). However, the same property does not hold when we consider the reach condition. It is not always possible to define the so-called $r$-hull of $S$, that is, the minimal set containing $S$ and having reach $\geq r$, see Colesanti and Manselli \cite{colesanti:10}. In fact, when $S$ admits such a minimal set it happens to coincide with $C_r(S)$, see Corollary 4.7 by Colesanti and Manselli \cite{colesanti:10}. This result provides an indirect proof of Proposition \ref{reachrconvex} below, that states that every set with reach $r$ is also $r$-convex.
Proposition  \ref{rconvexroll} establishes the relation between $r$-convexity and the rolling condition.

\begin{proposition}\label{reachrconvex} Let $S\subset
\mathbb{R}^d$ be a compact set with $\textnormal{reach}(S)\geq
r>0$. Then $S$ is $r$-convex.
\end{proposition}

\begin{proof} Let $x\in S^c$. If $\delta_S(x)\geq r$,
then $\mathring{B}(x,r)\cap S=\emptyset$. If $\delta_S(x)< r$, let
$\xi_S(x)$ be the unique point in $S$ such that
$\delta_S(x)=\left\|x-\xi_S(x)\right\|$. Then
$$
\eta_x=\frac{x-\xi_S(x)}{\left\|x-\xi_S(x)\right\|}\in
\textnormal{Nor}(S,\xi_S(x)),
$$
where $\textnormal{Nor}(S,\xi_S(x))$ is the set of all normal
vectors of $S$ at $\xi_S(x)$. Define
$y_\lambda=\xi_S(x)+\lambda\eta_x$ with $0<\lambda\leq r$ and take $\lambda\in(0,r)$. Now,
$\textnormal{reach}(S,\xi_S(x))\geq r>\lambda$ and by part (12) of Theorem 4.8 in Federer \cite{federer:59}, we get that $\eta_x=\lambda^{-1}
v$ where $\left\|v\right\|=\lambda$ and
$\xi_S(\xi_S(x)+v)=\xi_S(x)$. That is,
$\xi_S(y_\lambda)=\xi_S(x)$. It is
clear that $x\in\mathring{B}(y_r,r)$. Moreover,
$\mathring{B}(y_r,r)\cap S=\emptyset$. Indeed,  if
$z\in\mathring{B}(y_r,r)$, then
$z\in\mathring{B}(y_\lambda,\lambda)$ for some
$y_\lambda=\xi_S(x)+\lambda\eta_x$ with $0<\lambda<r$ and
$\mathring{B}(y_\lambda,\lambda)\cap S=\emptyset$ as a consequence
of $\xi_S(y_\lambda)=\xi_S(x)$ and $\mbox{reach}(S,\xi_S(x))\geq r>\lambda$.
\end{proof}

\begin{figure}[htb]
\begin{center}
\setlength{\unitlength}{1mm}
\begin{picture}(100,45)
\pscustom[linewidth=.5pt,fillstyle=solid, fillcolor=gray]{
\psarc(2,1.5){1}{90}{270}
\psline(2,.5)(0,.5)
\psline(0,.5)(0,1.5)
\psline(0,2.5)(2,2.5)
}

\pscustom[linewidth=.5pt,fillstyle=solid, fillcolor=gray]{
\psarc(3,1.5){1}{270}{90}
\psline(3,2.5)(5,2.5)
\psline(5,2.5)(5,0.5)
\psline(5,.5)(3,.5)
}

\psline[linewidth=.25pt,linestyle=dashed](1,1.5)(2,1.5)
\put(15,16){$r$}
\put(23,37){(a)}

\pscustom[linewidth=.5pt,fillstyle=solid, fillcolor=gray]{
\psarc(8,2){1}{180}{270}
\psarc(8,0){1}{90}{180}
\psline(7,0)(7,2)
}

\pscustom[linewidth=.5pt,fillstyle=solid, fillcolor=gray]{
\psline(9,0)(9,2)
\psarc(8,2){1}{0}{90}
\psline(8,3)(10,3)
\psarc(10,2){1}{90}{180}
}

\psline[linewidth=.5pt](9,0)(9,2)
\psline[linewidth=.25pt,linestyle=dashed](7,2)(8,2)
\put(75,21){$r$}
\put(80,37){(b)}

\end{picture}
\caption{(a) $r{\mbox{-convex}}\nRightarrow r{\mbox{-reach}}$. (b) $r{\mbox{-rolling}}\nRightarrow r{\mbox{-convex}}$.}
\label{fig:converse}
\end{center}
\end{figure}
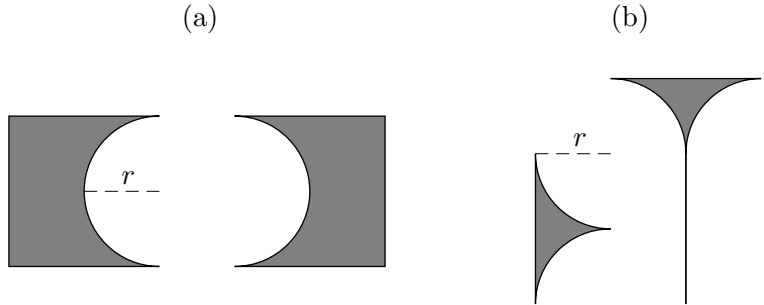

\begin{remark}\label{contract}
\it Borsuk's conjecture on local contractibility of $r$-convex sets\rm. The converse of Proposition \ref{reachrconvex} is not true in general, see Figure \ref{fig:converse} (a) for an example of a $r$-convex set that has not reach $r$. Even so, we prove in Theorem \ref{SSS} that if $S$ is a compact $r$-convex support in $\mathbb{R}^2$ fulfilling
a mild regularity condition (which we call \sl interior local connectivity\rm; see Section \ref{sec:appl} for details)
then  $S$ has positive reach, though not necessarily $r$. If we do not assume any additional regularity condition on $S$, then the conclusion of positive reach does not seem that simple to get. This is closely related to an unsolved conjecture by K. Borsuk (see Perkal \cite{perkal:56} and Mani-Levitska \cite{mani:93}): \sl Is an $r$-convex set locally contractible?\/ \rm   Note that proving that a compact $r$-convex set has positive reach would give a positive answer to Borsuk's conjecture since, according to  Remark 4.15 in Federer \cite{federer:59}, any set with positive reach is locally contractible. Recall that a topological space is said to be \sl contractible\/ \rm if it is homotopy equivalent to one point. In intuitive terms this means that the space can be continuously shrunk to one point. The space is called \sl locally contractible\/ \rm if every point has a local base of contractible neighborhoods.
\end{remark}

\begin{proposition}\label{rconvexroll} Let $S\subset
\mathbb{R}^d$ be a compact $r$-convex set for some $r>0$. Then $S$ fulfills the $r$-rolling condition.
\end{proposition}

\begin{proof} Let $x\in\partial S$. Since $x$ is a limit point of the set $S^c$, there is a sequence of points $\{x_n\}$, where $x_n\in S^c$, that converges to $x$. By the $r$-convexity, for each $n$ there exists $y_n\in\mathring{B}(x_n,r)$ such that $\mathring{B}(y_n,r)\cap S=\emptyset$. Since $\{y_n\}$ is bounded it contains a convergent subsequence which we denote by $\{y_n\}$ again.  Then $y_{n}\to y$ and it is not difficult to prove that $\mathring{B}(y,r)\cap S=\emptyset$ and $\left\|y-x\right\|\leq r$.  Now, $x\in S$ since $S$ is closed and, therefore, $x\in\partial {B}(y,r)$, which concludes the proof.
\end{proof}

The converse implication is not true in general, see Figure \ref{fig:converse} (b) for an example of a set fulfilling the $r$-rolling condition but not $r$-convex.

\section{Boundary convergence and full convergence in sequences of sets}\label{sub:seq}
As mentioned in the introduction, the focus in this paper  is on the reconstruction (in the statistical sense) of an unknown support $S$ from a sequence of estimators $\{S_n\}$ based on sample information. In many practical instances, including image analysis, the most important aspect of the target set is the boundary $\partial S$. However, it is clear that, even in very simple cases, the Hausdorff convergence $d_H(S_n,S)\to 0$ does not entail the
boundary convergence $d_H(\partial S_n,\partial S)\to 0$. See e.g. Ba\'{i}llo and Cuevas \cite{baillo:01}, Cuevas and Rodr\'{i}guez-Casal \cite{cuevas:04} and Rodr\'{i}guez-Casal \cite{rodriguez:07} for some results on boundary estimation. We introduce the following notion of convergence.

\begin{definition}\label{full}
Let $\{S_n\}$ be a sequence of compact non-empty sets in ${\mathbb R}^d$. Let $S\subset {\mathbb R}^d$ be a compact non-empty set. We say that we have full convergence of $\{S_n\}$ to $S$ if $d_H(S_n,S)\to 0$ and $d_H(\partial S_n,\partial S)\to 0$.
\end{definition}

The following result shows that, under very general conditions, the Hausdorff convergence of the sets and their boundaries implies also the convergence with respect to the distance in measure. This accounts for the term ``full convergence''.

\begin{theorem} \label{dHdmu}
Let $\{S_n\}$ be a sequence of compact non-empty sets in ${\mathbb R}^d$ endowed with a Borel measure $\nu$, with $\nu(C)<\infty$ for any compact $C$. Let $S$ be a compact non-empty set such that $\nu(\partial S)=0$, $d_H(S_n,S)\rightarrow 0$ and $d_H(\partial S_n,\partial S)\rightarrow 0$. Then $d_\nu(S_n, S)\rightarrow 0$.
\end{theorem}

\begin{proof} Take $\epsilon>0$. We first prove that, for $n$ large enough,
\begin{equation}
\{x\in S:\ \delta_{\partial S}(x)>2\epsilon\}\subset S_n\label{e1}
\end{equation}
and
\begin{equation}
\{x\in S_n:\ \delta_{\partial S_n}(x)>2\epsilon\}\subset S\label{e2}.
\end{equation}
To see (\ref{e1}) take $x\in S$ such that $\delta_{\partial S}(x)>2\epsilon$ and $n$ large enough such that
 $S\subset B(S_n,\epsilon)$ and $d_H(\partial S_n,\partial S)<\epsilon$. If $x\notin S_n$, then $x\in B(S_n,\epsilon)\setminus S_n$. Therefore $\delta_{\partial S_n}(x)\leq \epsilon$ so that
 $$
 \delta_{\partial S}(x)=d_H(\{x\},\partial S)\leq d_H(\{x\},\partial S_n)+d_H(\partial S,\partial S_n)<2\epsilon,
 $$
 which contradicts $\delta_{\partial S}(x)>2\epsilon$.

 The proof of (\ref{e2}) follows along the same lines.

 Finally, note that as a consequence of (\ref{e1}) and (\ref{e2}),
$$
S_n\Delta S\subset \{x\in S:\ \delta_{\partial S}(x)\leq 2\epsilon\}\cup \{x\in S_n:\ \delta_{\partial S_n}(x)\leq 2\epsilon\},
$$
which for large enough $n$ is a subset of  $\{x\in S: \ \delta_{\partial S}(x)\leq 3\epsilon\}$, that decreases to $\partial S$ as $\epsilon\downarrow 0$. Since $\nu$ is finite on bounded sets and $\nu(\partial S)=0$ this entails $\limsup \nu(S_n\Delta S)=0$ which concludes the proof.
\end{proof}

We next show that the above result applies to the important class of sets fulfilling the $r$-rolling condition. Therefore, as shown in Propositions \ref{reachrconvex} and \ref{rconvexroll}, it also applies to the class of $r$-convex sets and to that of sets with positive reach.

\begin{theorem}\label{partialSn}
Let $\{S_n\}$ be a sequence of compact non-empty sets in $\mathbb{R}^d$ satisfying the $r$-rolling condition. Let $S$ be a compact non-empty set such that $d_H(S_n,S)\rightarrow 0$. Then,

(a)  $d_H(\partial S_n,\partial S)\rightarrow 0$.

(b) If $\textnormal{reach}(S_n)\geq r>0$, then we also have $d_\nu(S_n,S)\rightarrow 0$ for any Borel measure $\nu$ (finite on compacts) absolutely continuous with respect to the Lebesgue measure $\mu$.
\end{theorem}

\begin{proof}(a) Assume that the result is not true. Then,
\begin{itemize}
\item[(i)]{There exists $\epsilon>0$ such that for infinitely many  $n\in\mathbb{N}$ there exists $x_n\in\partial S$ with $\delta_{\partial S_n}(x_n)>\epsilon$ or}
\item[(ii)]{There exists $\epsilon>0$ such that for infinitely many $n\in\mathbb{N}$ there exists $x_n\in\partial S_n$ with $\delta_{\partial S}(x_n)>\epsilon$.}
\end{itemize}

First, assume that (i) is satisfied. Since $S$ is compact, there exists a convergent subsequence of $\{x_n\}$ which we will denote again as the original sequence. Let $x\in \partial S$ be the limit of $\{x_n\}$. The Hausdorff convergence of $S_n$ to $S$ implies that $\delta_{S_n}(x)\leq\epsilon/2$ for infinitely many $n$. Furthermore, for large enough $n$,
\[\delta_{\partial S_n}(x)\geq \delta_{\partial S_n}(x_n)-\left\|x-x_n\right\|> \epsilon/2.\]
Now, $\delta_{S_n}(x)\leq\epsilon/2$ together with $\delta_{\partial S_n}(x)>\epsilon/2$ yields $x\in\textnormal{int}(S_n)$ for infinitely many $n$ and $\mathring{B}(x,\epsilon/2)\cap S_n^c=\emptyset$. But, since $x\in\partial S$ and $S$ is closed, we can consider $y\in S^c$ such that $\left\|x-y\right\|<\epsilon/2$ and $\delta_S(y)>0$. Again by the Hausdorff convergence of $S_n$ to $S$ we get that $y\in S_n^c$ for infinitely many $n$ which yields a contradiction.

Assume now that (ii) is satisfied. We can assume $\epsilon<r$. The
Hausdorff convergence of $S_n$ to $S$ implies that
$\delta_{S}(x_n)\leq\epsilon$ for infinitely many $n$. This, together
with $\delta_{\partial S}(x_n)>\epsilon$, yields
$x_n\in\textnormal{int}(S)$ for infinitely many $n$.
Let $0<\lambda<\epsilon/2<r$.
Since $x_n\in\partial S_n$ and the sets $S_n$ satisfy the $r$-rolling property, there exists
for each $n\in\mathbb{N}$ a ball ${B}(c_n,\lambda)$ such that $x_n\in\partial
{B}(c_n,\lambda)$ and $\mathring{B}(c_n,\lambda)\cap
S_n=\emptyset$. Again, let us denote by $\{x_{n}\}$ a convergent subsequence of
$\{x_n\}$. Then $x_{n}\to x$, with $x\in
S$ and $\delta_{\partial S}(x)\geq\epsilon/2$, that is,
$B(x,\epsilon/2)\subset S$. Now, let $\{c_{n}\}$ be a
convergent subsequence of $\{c_{n}\}$. Then
$c_{n}\to c$, with $c\in S_n^c$  and
$\delta_{S_{n}}(c)>\lambda/2$ for infinitely many $n$. By the Hausdorff
convergence of $S_n$ to $S$, we get that $c\in S^c$ which yields a
contradiction with $B(x,\epsilon/2)\subset S$.

(b) By Theorem 4.13 in Federer \cite{federer:59}, $\textnormal{reach}(S)\geq r$ and, in particular, $\mu(\partial S)=0$, see Remark \ref{remmc}. The result is now straightforward from (a) and Theorem \ref{dHdmu}.
\end{proof}

\begin{remark}\label{compact}
If ${\mathcal A}$ is a class of sets, closed with respect to the Hausdorff topology, and the $d_H$-convergence in this class implies the $d_H$-convergence of the respective boundaries, then Theorem \ref{dHdmu} entails that this class is also closed with respect to $d_\nu$
whenever $\nu(\partial A)=0$ for all $A\in{\mathcal A}$. So, the $d_H$-compactness of ${\mathcal A}$ implies the $d_\nu$-compactness of this class.

On the other hand, Theorem \ref{partialSn} shows that this in fact applies to the class of sets with \textnormal{reach} $\geq r$, which is $d_H$-closed from Theorem 4.13 and Remark 4.14 in Federer \cite{federer:59}. Thus, the class of subsets with \textnormal{reach} $\geq r$ of a compact set is $d_\mu$-compact.
This will be useful below (subsection \ref{EM}) in order to apply the results in Polonik \cite{polonik:95} for classes of sets defined in terms of reach properties.
\end{remark}

\section{$P$-uniformity: Billingsley-Tops{\o}e theory and its application to classes of sets with positive reach}\label{sec:bill}

Let $X_1, X_2,\ldots$ be a sequence of independent and identically distributed random elements on a probability space (${\Omega}$, $\mathcal{F}$, $\mathbb{P}$) with values in a measurable space ($E$, $\mathcal{B}$), where $E$ is a metric space and $\mathcal{B}$ stands for the Borel $\sigma$-algebra on $E$. Denote by $P$
the probability distribution of $X_1$ on $\mathcal{B}$ and let  ${\mathbb P}_n$ be the empirical probability measure associated with $X_1,\ldots,X_n$.

The almost sure pointwise convergence on $\mathcal{B}$ of ${\mathbb P}_n$ to $P$ is ensured by the
strong law of large numbers. Moreover, for appropriate classes ${\mathcal A}\subset {\mathcal B}$ the uniform convergence
\begin{equation}
\sup_{A\in{\mathcal A}}|{\mathbb P}_n(A)-P(A)|\to 0,\ \mbox{a.s.}\label{GC}
\end{equation}
also holds.

A class ${\mathcal A}$ of sets fulfilling (\ref{GC}) is called a {\slshape{Glivenko-Cantelli class}} (GC-class). They are named after the classical Glivenko-Cantelli theorem which establishes the result (\ref{GC}) for the case where $E={\mathbb R}$ and ${\mathcal A}$ is the class of closed half lines $A=(-\infty,x]$, $x\in{\mathbb R}$. The study of uniform results of type (\ref{GC}) is a classical topic in statistics. A well-known reference is Pollard \cite{pollard:84}. A useful summary, targeted to the most usual applications  in statistics, can be found in Chapter 19 of van der Vaart \cite{vaart:98}.

A popular methodology to obtain Glivenko-Cantelli classes is based on the use of the well-known Vapnik-Cervonenkis inequality (see, e.g. Devroye et al. \cite{devroye:96}) which relies on combinatorial tools. This inequality provides upper bounds for ${\mathbb P}\{\sup_{A\in{\mathcal A}}|{\mathbb P}_n(A)-P(A)|>\epsilon\}$ which depend on the so\--called \sl shatter coefficients \/ \rm and \sl VC-dimension\/ \rm  of the class ${\mathcal A}$. However, this approach is not useful in those cases where the VC-dimension of ${\mathcal A}$ is infinite. This is the case where ${\mathcal A}$ is the class of closed convex sets in ${\mathbb R}^d$ and therefore also for the class of closed $r$-convex sets.

Nevertheless, it can be shown that the family of all convex sets in ${\mathbb R}^d$ is a GC-class. This can be done following an alternative approach (maybe less popular than the VC-method), due to
 Billingsley  and Tops{\o}e \cite{billingsley:67}; see also Bickel and Millar \cite{bickel:92}. This methodology is rather based on geometrical and topological ideas and therefore turns out to be more suitable for set estimation purposes. The basic ideas of Billingsley-Tops{\o}e approach can be summarized as follows.

Let $P_n$ and $P$ be probability
measures on $\mathcal{B}$. A set $A$ is called a {\slshape{$P$-continuity
set}} if $P(\partial A)=0$. The sequence $P_n$ is said to converge
weakly to $P$ if
\begin{equation*}\label{Pweak}
P_n\left(A\right)\to P\left(A\right)
\end{equation*}
for each $P$-continuity set $A\subset\mathcal{B}$. A subclass $\mathcal{A}\subset\mathcal{B}$ is said to be a {\slshape{$P$-continuity class}} if every set in $\mathcal{A}$ is a $P$-continuity set. A subclass $\mathcal{A}\subset\mathcal{B}$ is said to be a {\slshape{$P$-uniformity class}} if
\begin{equation*}\label{Punif}
\sup_{A\in\mathcal{A}}\left|P_n\left(A\right)-P\left(A\right)\right|\to 0
\end{equation*}
holds for every sequence $P_n$ that converges weakly to $P$. Note that the $P$-uniformity concept is not established just for sequences of empirical distributions but in general for any sequence of probability measures converging to $P$. Billingsley  and Tops{\o}e \cite{billingsley:67} derived several results establishing conditions for a class $\mathcal{A}$ to be a $P$-uniformity class. As a consequence, some useful criteria for obtaining GC-classes immediately follow.

The following theorem provides three sufficient conditions to ensure that a class ${\mathcal A}$ is a $P$-uniformity class.

\begin{theorem}\label{Bill}(Billingsley and Tops{\o}e (1967, Th. 4))
If $E$ is locally connected and if $\mathcal{A}$ is a $P$-continuity class of subsets of $E$, then each of the following three conditions is sufficient for $\mathcal{A}$ to be a $P$-uniformity class.
\begin{itemize}
\item[(i)]{The class $\partial\mathcal{A}=\{\partial A:\ A\in\mathcal{A}\}$ is a compact subset of the space $\mathcal{M}$ of non-empty closed bounded subsets of $E$.}
\item[(ii)]{There exists a sequence $\{C_n\}$ of bounded sets with $P(\textnormal{int}(C_n))\to 1$ and such that, for each $n$, the class
$\partial(C_n\cap\mathcal{A})=\{\partial(C_n\cap{A}):\ A\in\mathcal{A}\}$
is a compact subset of $\mathcal{M}$.}
\item[(iii)]{There exists a sequence $\{C_n\}$ of closed, bounded sets with $P(\textnormal{int}(C_n))\to 1$ and such that, for each $n$, the class
$C_n\cap\partial\mathcal{A}=\{C_n\cap\partial{A}:\ A\in\mathcal{A}\}$
is a compact subset of $\mathcal{M}$. }
\end{itemize}
\end{theorem}
\rm

\

Theorem \ref{Bill} can be used to prove that the class of all convex sets in ${\mathbb R}^d$ is a $P$-uniformity class, and in particular, a Glivenko-Cantelli class. The following theorem provides a partial extension of this property: we show that, under an additional boundedness assumption, convexity can be replaced with the broader condition of having a given positive reach. Again, the basic tool is Theorem \ref{Bill}.

\begin{theorem}\label{reachcompact}
Let $K$ be a compact non-empty subset of $\mathbb{R}^d$ and $\mathcal{A}=\{A\subset
K:\ A\neq\emptyset,\ A \textit{ is closed, and } \textnormal{reach}(A)\geq r\}$, for $r>0$. Then the class
$\partial\mathcal{A}=\{\partial A:\ A\in\mathcal{A}\}$ is a compact subset of $\mathcal{M}$. Moreover, $\mathcal{A}$ is a $P$-uniformity class for every probability measure $P$ such that $P$ is absolutely continuous with respect to the Lebesgue measure $\mu$.
\end{theorem}

\begin{proof} Let $\{\partial A_n\}$ be a convergent sequence of sets in
$\partial\mathcal{A}$. By
Federer's closeness theorem for sets of positive reach (Theorem
4.13 in Federer \cite{federer:59}) it follows that the class $\mathcal{A}$
is compact with respect to the Hausdorff metric and, therefore,
$A_n$ has a convergent subsequence whose limit is a set $A\in\mathcal{A}$. In view of Propositions \ref{reachrconvex} and \ref{rconvexroll}, we can apply Theorem \ref{partialSn} to the class $\mathcal{A}$, and this yields $\partial{A_n}\to \partial {A}$ in $d_H$.
Now, from Theorem \ref{Bill} (i), $\mathcal{A}$ is a $P$-uniformity class.
\end{proof}

\section{Applications}\label{sec:appl}

\subsection{Estimation of $r$ convex sets and their boundary lengths}\label{sub:reach}

\noindent \it Estimation of $r$-convex supports\rm

\begin{figure}[h]
\begin{center}
\setlength{\unitlength}{1mm}
\begin{picture}(60,65)
\put(0,0){\scalebox{0.3}{\includegraphics{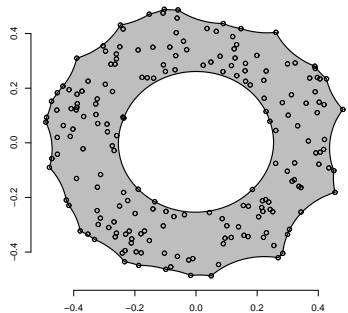}}}
\end{picture}
\caption{In gray, $r$-convex hull of a sample of size $n=200$ from a uniform distribution on the set $S=B(0,0.5)\setminus\mathring{B}(0,0.25)$. The value of $r$ is 0.25.}
\label{fig:rh}
\end{center}
\end{figure}

As indicated in Section \ref{rolling}, $r$-convexity is a natural extension of the notion of convexity. From the point of view of set estimation, $r$-convexity is particularly attractive as the estimation of an $r$-convex compact support $S$ from a random sample $X_1,\ldots,X_n$ drawn on $S$, can be handled very much in the same way as the case where $S$ is convex.  In this classical situation (which has been extensively considered in the literature, see, e.g., D{\"u}mbgen and Walther \cite{dumbgen:96}) the natural estimator of $S$ is the convex hull of the sample. In an analogous way, if $S$ is assumed to be $r$-convex, the obvious estimator is the $r$-convex hull of the sample, which we will denote by $S_n$. Recall from (\ref{rhull}) that
\begin{equation}\label{rhulln}
  S_n=\bigcap_{\mathring{B}(y,r)\cap
    \{X_1,\ldots,X_n\}=\emptyset}\mathring{B}(y,r)^c.
\end{equation}
This estimator can be explicitly calculated in a computationally efficient way (at least in the two-dimensional case) through the R-package \tt alphahull\rm; see Pateiro-L\'{o}pez and Rodr\'{i}guez-Casal \cite{pateiro:10}. In Figure \ref{fig:rh} we show an example of the $r$-convex hull estimator.
Note that the boundary of the $r$-convex hull estimator is formed by arcs of balls of radius $r$ (besides
possible isolated sample points). The arcs are determined by the intersections of some
of the empty balls that define the complement of the $r$-convex hull, see Equation (\ref{rhulln}).

The $r$-convex hull estimator was first considered, from a statistical point of view, in Walther \cite{walther:97}. Rodr\'{i}guez-Casal \cite{rodriguez:07} studied its properties as an estimator of $S$ and $\partial S$, providing (under some rolling-type assumptions for $S$) convergence rates for $d_H(S_n,S)$, $d_H(\partial S_n,\partial S)$ and $d_\mu(S_n,S)$ which essentially coincide with those given by D{\"u}mbgen and Walther \cite{dumbgen:96} for the convex case.

\

\noindent \it Estimation of the boundary length\rm

We are interested in estimating the boundary length of $S$, $L(S)$, as defined by the outer Minkowski content given in (\ref{os}). Recall that (see Remark \ref{remmc}) if $\textnormal{reach}(S)>0$, then the outer Minkowski content is finite.

In Pateiro-L\'{o}pez and Rodr\'{i}guez-Casal \cite{pateiro:08} the estimation of the usual (two-sided) Minkowski content $L_0(S)$ given in (\ref{mc}) is considered under double smoothness assumptions of rolling-type on $S$ (from inside and outside $\partial S$) which in fact are stronger than $r$-convexity. Under these assumptions, these authors improve the convergence rates for the Minkowski content $L_0(S)$ obtained in Cuevas et al. \cite{cuevas:07}. Another recent contribution to the problem of estimating boundary measures is due to
Jim\'{e}nez and Yukich \cite{jimenez:10}. In all these cases the estimators are based on a sample model which requires random observations inside and outside $S$ in such a way that for each observation one is able to decide (with no error) whether or not  it belongs to $S$.

Theorem \ref{SSS} provides further insights on the approach of these papers in the sense that, under minimal assumptions, gives a fully consistent estimator $S_n$ of an $r$-convex set $S\subset{\mathbb R}^2$ and a plug-in consistent estimator $L(S_n)$ of $L(S)$ based on a unique inside sample.
The performance of $L(S_n)$, in terms of bias and variance is analyzed in the appendix.

 Before stating Theorem \ref{SSS} we need some preliminaries.

\begin{definition}
We will say that a set $S$ fulfills the property of interior local connectivity (ILC) if there exists $\alpha_0>0$ such that for all $\alpha\leq\alpha_0$ and for all $x\in S$, $\textnormal{int}(B(x,\alpha)\cap S)$ is a non-empty connected set. \rm
\end{definition}

\begin{figure}[htb]
\begin{center}
\setlength{\unitlength}{1mm}
\begin{picture}(50,25)
\pscustom[linewidth=.5pt,fillstyle=solid, fillcolor=gray]{
\psarc(2.5,2){2}{180}{360}
\psline(4.5,2)(4.5,0)
\psline(4.5,0)(0.5,0)
\psline(0.5,0)(0.5,2)
}
\put(25,2){$x$}
\psdot(2.5,0)
\put(6.5,4){$S$}
\end{picture}
\caption{The set $S$ in gray does not fulfill the ILC property at the central point $x$.}
\label{fig:ILC}
\end{center}
\end{figure}
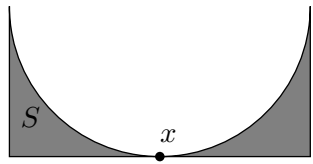

\

\begin{figure}[htb]
\begin{center}
\setlength{\unitlength}{1mm}
\begin{picture}(10,45)

\psarc[linestyle=dashed,linewidth=.5pt](0,1){1}{0}{360}
\psarc[linestyle=dashed,linewidth=.5pt](0.75,2.3){1}{0}{360}
\pscustom[linewidth=.5pt,fillstyle=solid, fillcolor=gray]{
\psarc(0,1){1}{101.4}{150}
\psarc(0.75,2.3){1}{150}{198.59}
}
\psdot[dotsize=0.1](0.75,2.3)
\psline[linestyle=dashed,linewidth=.5pt](0.75,2.3)(-0.25,2.3)
\put(2.5,24){$r$}

\end{picture}
\caption{In gray, an $r$-circular triangle.}
\label{fig:rtri}
\end{center}
\end{figure}
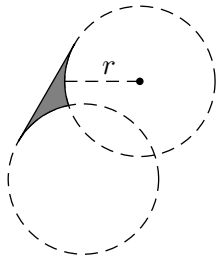

Let us define an \sl $r$-circular triangle\/ \rm  as a compact plane figure limited by three sides: two of them are arcs of intersecting circumferences of radius $r$; the third side is a linear segment which is tangent to both circumference arcs. In other words, an $r$-circular triangle is an isosceles triangle with a linear segment as a basis and two $r$-circumference segments as the other sides. See Figure \ref{fig:rtri}.

\begin{theorem}\label{SSS}
\sl Let $S\subset{\mathbb R}^2$ be a compact $r$-convex set with $\mu(S)>0$ that fulfills the ILC property. Let $X_1,\ldots,X_n$ be a random sample drawn from a uniform distribution with support $S$. Denote by $S_n$ the $r$-convex hull of this sample, as defined in (\ref{rhulln}). Then,

(a) $S_n$ is a fully consistent estimator of $S$.

(b) $S$ has positive reach and $L(S_n)\to L(S)$,  a.s., $L(S)$ being the outer Minkowski content of $S$.\rm

\end{theorem}

\begin{proof} (a) Since $\{X_1,\ldots,X_n\}\subset S_n\subset S$ and $d_H(\{X_1,\ldots,X_n\},S)\to 0$, a.s. we also have $d_H(S_n,S)\to 0$, a.s. In order to get the full consistency we only need to prove
$d_H(\partial S_n,\partial S)\to 0$, a.s. but this follows directly from Proposition \ref{rconvexroll} and
Theorem \ref{partialSn}.

\

\noindent(b) Define $\tilde{S_n}=S_n\setminus I(S_n)$, where $I(S_n)$ is the set of isolated points of $S_n$. Thus, for large enough $n$, with probability one $\tilde{S_n}$ is a compact (not necessarily connected) set whose boundary is the union of a finite sequence of $r$-arcs, that is, circumference arcs of radius $r$.
We will say that a  boundary point of $\tilde{S_n}$ is  ``extreme'' if it is the intersection of two boundary $r$-arcs from different $r$-circumferences. The proof is based on the following two lemmas.

\begin{lemma}\label{lemareach}
\sl Under the conditions of Theorem \ref{SSS}, there exists $r_0>0$ such that $\textnormal{reach}(\tilde{S_n})\geq r_0$ a.s. for all $n$.
\end{lemma}

\begin{proof} If the sequence $\mbox{reach}(\tilde{S_n})$ were not bounded from below a.s. we would have, with positive probability, a sequence $\{z_n\}$ of points   with $z_n\notin \tilde{S_n}$ and another  sequence  $(x_n,y_n)$, where $x_n$ and $y_n$ are boundary points of $\tilde{S_n}$, with $x_n\neq y_n$ and such that
\begin{equation}
\Vert z_n-x_n\Vert=\Vert z_n-y_n\Vert=\inf_{x\in \tilde{S_n}}\Vert z_n-x\Vert=r_n\to 0.\label{cont}
\end{equation}
In principle, the ``projection points'' $x_n$ and $y_n$ could be extreme points (as defined above) or ``boundary inside points'', that is, points belonging to a boundary $r$-arc but different from the extremes of this arc. However, it is easily seen (Figure \ref{fig:arc}) that the latter possibility (i.e.,  $x_n$ or $y_n$ is a boundary inside point) leads to a contradiction: Indeed, note that if $x_n$ belongs to an $r$-arc in $\partial \tilde{S_n}$ with extremes $a$ and $b$ and $x_n\neq a$, $x_n\neq b$, then at a fixed distance $r_n<r$ there is only one point
in the complement of $\tilde{S_n}$ (necessarily equal to $z_n$) such that $x_n$ is the projection of $z_n$ on $\tilde{S_n}$. This would imply that the projection of $z_n$ is unique since the arc with extremes $a$ and $b$ is an arc of a circle of radius $r$ that does not intersect $\tilde{S_n}$.

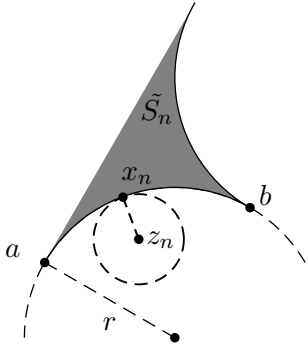
\begin{figure}[htb]
\begin{center}
\setlength{\unitlength}{1mm}
\begin{picture}(50,60)
\psarc[linestyle=dashed,linewidth=.5pt](2.46,1){2}{30}{180}
\psdot(2.46,1)
\psline[linestyle=dashed,linewidth=.5pt](2.46,1)(0.73,2)
\put(15,11){$r$}

\pscustom[linewidth=.5pt,fillstyle=solid, fillcolor=gray]{
\psarc(4.46,4.46){2}{150}{240}
\psarc(2.46,1){2}{60}{150}
}

\psdot(0.73,2)
\psdot(3.46,2.73)

\psdot(1.77,2.87)
\psdot(1.98,2.31)
\psline[linestyle=dashed](1.77,2.87)(1.98,2.31)
\psarc[linestyle=dashed, linewidth=0.025](1.98,2.31){0.6}{0}{360}

\put(35.8,27.7){$b$}
\put(2,20.8){$a$}
\put(17.5,31){$x_n$}
\put(20.8,22.5){$z_n$}

\put(20,39){$ \tilde{S_n}$}

\end{picture}
\caption{If $z_n$ has two projections $x_n$ and $y_n$ onto $\tilde{S_n}$, neither $x_n$ nor $y_n$ can be boundary inside points.}
\label{fig:arc}
\end{center}
\end{figure}

Thus, $x_n$ and $y_n$ must be extreme points and our proof reduces to see that we cannot have (\ref{cont}) with a positive probability. So, let us assume that (\ref{cont}) were true with a positive probability (let us denote by $A$ the corresponding event).

Let $T_n^1$ and $T_n^2$ be two circular triangles, with vertices $x_n$ and $y_n$, respectively, determined by the $r$-arcs in $\partial \tilde{S_n}$ whose intersections are $x_n$ and $y_n$, see Figure \ref{fig:tn1}. Note that these triangles are not necessarily included in $\tilde{S_n}$; just a
portion of each triangle close to the vertex is in general included in  $\tilde{S_n}$.

First, let us prove that the  ``heights'' (i.e. the distances from the vertices to the bases) of $T_n^1$ and $T_n^2$ must necessarily be bounded from below a.s. on $A$. To see this, recall that each $x_n$ is the intersection of the boundary arcs of two balls of radius $r$, $B(c_{1n},r)$ and $B(c_{2n},r)$, whose interiors do not intersect $\tilde{S_n}$. By construction, $\mbox{height}(T_n^1)\to 0$ for some subsequence $T_n^1$ implies that the centers $(c_{1n},c_{2n})$ must fulfill $c_1=\lim c_{1n}=\lim c_{2n}=c_2$. This means that the $r$-circumferences providing the boundary of $\tilde S_n$ at both sides of $x_n$ tend to coincide as $n$ tends to infinity. We will prove that $c_1=c_2$ leads to a contradiction. Indeed, if $c_1=c_2$ we would have that, for large enough $n$, the centers $c_{1n}$ and $c_{2n}$ would be very close and the boundary arcs, $\partial B(c_{1n},r)$ and $\partial B(c_{2n},r)$, would intersect to each other not only at $x_n$ but also at another point $x_n^\ast$ such that $\left\|x_n-x_n^\ast\right\|>\gamma>0$. We conclude $B(z_n,r_n)\setminus\{x_n\}\subset \textnormal{int}(B(c_{1n},r)\cup B(c_{2n},r))\subset \tilde{S_n}^c$, for large enough $n$ (to see this note that the boundary of $B(c_{1n},r)\cup B(c_{2n},r)$ near $x_n$ coincides with the triangular sides of $T_n^1$). Now, from (\ref{cont}) we have obtained a contradiction $y_n\in\partial B(z_n,r_n)\cap \tilde S_n$ and $y_n\neq x_n$.

\begin{figure}
\begin{center}
\setlength{\unitlength}{1mm}
\begin{picture}(80,60)
\put(-10,-15){\scalebox{0.45}{\includegraphics{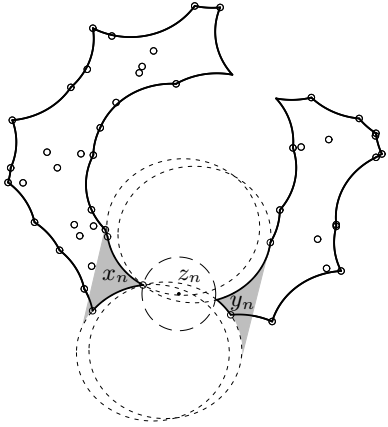}}}
\psdot[dotsize=0.05](3.52,1.41)
\psarc[linestyle=dashed, linewidth=0.01](3.52,1.41){0.49}{0}{360}
\put(35,16){\footnotesize{$z_n$}}
\put(25,16){\footnotesize{$x_n$}}
\put(42,12.5){\footnotesize{$y_n$}}
\end{picture}
\caption{In gray, $T_n^1$ and $T_n^2$.}
\label{fig:tn1}
\end{center}
\end{figure}

 As the space of compact sets endowed with the Hausdorff metric is locally compact, there exist a.s. convergent subsequences of $T_n^1$ and $T_n^2$, $\{x_n\}$ and $\{y_n\}$ which we will denote again as the original sequences.

 Let $T^1$ and $T^2$ be the a.s. Hausdorff limits of $T_n^1$ and $T_n^2$, respectively. Since the heights of $T_n^1$ and $T_n^2$ are bounded from below, $T^1$ and $T^2$ must be also non-degenerate circular triangles. Denote $x=\lim x_n=\lim y_n=\lim z_n$, a.s. Note that the fact that $z_n\notin S_n$ ensures that $x\in \partial S$. Then, we have two possibilities,

 (i) if $T^1\cap T^2=\{x\}$ we would get a contradiction: to see this note that $T=T^1\cup T^2$ is just the union of two $r$-circular triangles which are disjoint except for the common vertex $\{x\}$. Thus, for $\epsilon$ small enough $B(x,\epsilon)\cap S$ is included in $T$. This contradicts the ILC of $S$ at $x$. See Figure \ref{fig:limit} (a).

\begin{figure}[htb]
\begin{center}
\setlength{\unitlength}{1mm}
\begin{picture}(85,55)

\pscustom[linewidth=.5pt,fillstyle=solid, fillcolor=gray]{
\psarc(2,1){1}{110}{200}
\psline(1.07,0.66)(0.38,2.53)
\psarc(1.32,2.87){1}{200}{290}
}

\pscustom[linewidth=.5pt,fillstyle=solid, fillcolor=gray]{
\psarc(1.32,2.87){1}{290}{380}
\psline(2.25,3.22)(2.93,1.34)
\psarc(2,1){1}{20}{110}
}

\put(16,21){$x$}
\psdot[dotsize=0.07](1.68,1.95)

\pscustom[linewidth=.5pt,fillstyle=solid, fillcolor=gray]{
\psarc(7,1){1}{110}{200}
\psline(6.07,0.66)(5.38,2.53)
\psarc(6.32,2.87){1}{200}{290}
}

\pscustom[linewidth=.5pt,fillstyle=solid, fillcolor=gray]{
\psarc( 6.72, 1.67){1}{90}{180}
\psline(5.73 ,1.67)(  5.72 ,3.66)
\psarc(6.72, 3.67){1}{180}{270}
}

\pscustom[linewidth=0.2pt, fillstyle=hlines, hatchwidth=0.2pt]{
\psarc( 6.72, 1.67){1}{0}{360}
}

\put(68,27){\footnotesize{$x_n$}}
\put(68,18){\footnotesize{$y_n$}}
\put(65,10){\footnotesize{$\tilde{S_n}^c$}}
\psdot[dotsize=0.07](6.7,2.67)
\psdot[dotsize=0.07](6.65,1.935)

\put(57.5,26){\footnotesize{$T_n^1$}}
\put(57.5,15.8){\footnotesize{$T_n^2$}}
\psdot[dotsize=0.07](1.68,1.95)
\put(12,45){(a)}
\put(60,45){(b)}

\end{picture}
\caption{(a) $T_1\cap T_2=\{x\}$ yields a contradiction with the ILC property. (b) $T_1\cap T_2\neq\{x\}$ yields a contradiction with the fact that the circles of radius $r$ defining the arcs of $T_n^1$ and $T_n^2$ cannot intersect $\tilde{S_n}$.}
\label{fig:limit}
\end{center}
\end{figure}
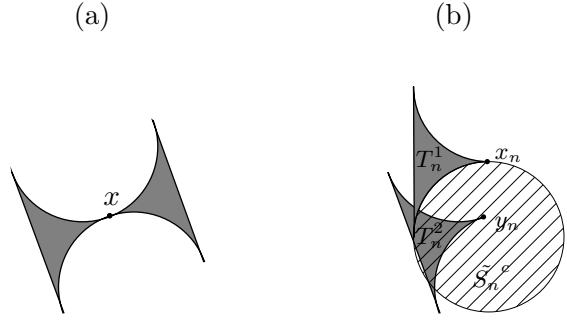

(ii) On the other hand, the possibility $T^1\cap T^2\neq \{x\}$ leads also to contradiction. As the vertices of $T_n^1$ and $T_n^2$ tend to the same point, we would necessarily have that, with probability one, $y_n$ (or $x_n$) belongs to the interior of one of the $r$-circles defining the arcs of $T_n^1$ (or $T_n^2$), see Figure \ref{fig:limit} (b). However, by construction, the $r$-circles defining the arcs of $T_n^1$ and $T_n^2$ do not intersect $\tilde{S_n}$.
This concludes the proof of Lemma \ref{lemareach}.
\end{proof}

\begin{lemma}\label{Isol1}
\sl Under the conditions of Theorem \ref{SSS},
\begin{equation}\label{hausIS}
d_H(\tilde{S_n},S)\to 0,\ \mbox{a.s.}
\end{equation}
\end{lemma}

\begin{proof} Assume that (\ref{hausIS}) is not true. Then, with positive probability we have either
\begin{itemize}
\item[(i)]{There exists $\epsilon>0$ such that for infinitely many $n\in\mathbb{N}$ there exists $x_n\in \tilde{S_n}$ with $\delta_{S}(x_n)>\epsilon$ or}
\item[(ii)]{There exists $\epsilon>0$ such that for infinitely many $n\in\mathbb{N}$ there exists $x_n\in S$ with {\mbox{$\delta_{\tilde{S_n}}(x_n)>\epsilon$}}.}
\end{itemize}

First, since $S$ is $r$-convex, $S_n\subset S$ with probability one. Therefore $\tilde{S_n}\subset S_n\subset B(S,\epsilon)$ for all $\epsilon>0$ and we cannot have (i).

With regard to (ii), let $\{x_n\}$ be a sequence of points in $S$. Since $S$ is compact, there exists a convergent subsequence which we will denote again $\{x_n\}$. Let $x\in S$ be the limit of $\{x_n\}$ and $y\in \textnormal{int}(S)$ such that $\left\|x-y\right\|<\epsilon/2.$ Note that the existence of such $y$ follows from the ILC property. Then, there exists $\epsilon_0>0$ ($\epsilon_0\ll r$) such that  $B(y,\epsilon_0) \subset S$. Using (a) in Theorem \ref{SSS} and a similar reasoning to that of (\ref{e1}) we get that, for $n$ large enough, $B(y,\epsilon_1) \subset S_n$ for some $\epsilon_1<\epsilon_0$. Finally, for $n$ large enough, 
\[\delta_{\tilde{S_n}}(x_n)\leq \delta_{B(y,\epsilon_1)}(x_n)\leq\delta_{B(y,\epsilon_1)}(x)+\left\|x_n-x\right\|\leq\left\|x-y\right\|+\left\|x_n-x\right\|<\epsilon.\]
That is, (ii) neither can be true. This concludes the proof of Lemma \ref{Isol1}.
\end{proof}

Now, in view of Lemmas \ref{lemareach} and \ref{Isol1}, the assumptions of Theorem 5.9 in Federer \cite{federer:59} are fulfilled (see also Remark 4.14 in that paper).  Theorem 5.9 in Federer \cite{federer:59} essentially establishes that $S$ has positive reach and the curvature measures are continuous with respect to $d_H$ (see Remark 5.10 in Federer \cite{federer:59}).
In particular we obtain that $\Phi_{d-1}(\tilde{S}_n, K)\to \Phi_{d-1}(S, K)$ a.s. for any closed ball $K$ such that $S\subset K$.  Using Remark 5.8 in Federer \cite{federer:59} and $\tilde{S}_n\subset S$ we get that $\Phi_{d-1}(\tilde{S}_n, K)=\Phi_{d-1}(\tilde{S}_n, K\cap\partial\tilde{S}_n)=\Phi_{d-1}(\tilde{S}_n, \partial \tilde{S}_n)$ and also $\Phi_{d-1}(S, K)=\Phi_{d-1}(S,\partial S)$. The proof of Theorem \ref{SSS} (b) concludes noting (see Remark \ref{remmc}) that $L(S_n)=L(\tilde{S}_n)=\Phi_{d-1}(\tilde{S}_n, \partial \tilde{S}_n)$ and $L(S)=\Phi_{d-1}(S, \partial S)$.

\end{proof}

\begin{remark}\label{Borsuk}
This result provides, using stochastic methods, a partial converse of Proposition \ref{reachrconvex}: we prove that $r$-convexity implies positive reach for ILC sets in ${\mathbb R}^2$. Thus we also get a partial answer to Borsuk's question (is  an $r$-convex set locally contractible?) mentioned in Remark \ref{contract}, since from Federer \cite{federer:59}, Remark 4.15, any set with positive reach is locally contractible.
\end{remark}

\subsection{Applications to the excess mass approach}\label{EM}

Typically, the results on uniform convergence, of Glivenko-Cantelli type, are useful in order to establish the
consistency of set estimators which are defined as maximizers of appropriate functionals of the empirical process.
An interesting example in set estimation is the \sl excess mass approach\rm, proposed by Hartigan \cite{hartigan:87}
and M\"uller and Sawitzki \cite{muller:91} and further developed by Polonik \cite{polonik:95} and Polonik and Wang \cite{polonik:05}, among others.

The basic ideas of this method can be simply described as follows: given $\lambda>0$, denote by $P$ an
absolutely continuous distribution with respect to the Lebesgue measure $\mu$. Let $f$ denote the
corresponding $\mu$-density. Define the \sl excess mass functional\/ \rm on the class ${\mathcal B}({\mathbb R}^d)$ of Borel sets,
\begin{equation*}
H_\lambda(A)=P(A)-\lambda \mu(A),\ \mbox{for\ } A\in {\mathcal B}({\mathbb R}^d).\label{H}
\end{equation*}
Since $H_\lambda(A)=\int_A(f-\lambda)d\mu$, it is obvious that $H_\lambda(A)$ is maximized by the level set $A=\{f\geq \lambda\}$.

This suggests a method to define an estimator of
$S(\lambda)=\{f\geq\lambda\}$ which can incorporate some shape
restrictions previously imposed on this set. Let
us assume that $S(\lambda)$ belongs to some given class of sets
${\mathcal A}$ (for example, the class of compact convex sets or
the class of compact $r$-convex sets in ${\mathbb R}^d$). As
$S(\lambda)$ is the maximizer on ${\mathcal A}$ of the unknown
functional $H_\lambda(A)$, we could define $S_n(\lambda)$ as the
maximizer on ${\mathcal A}$ of the \sl empirical excess mass
functional\rm
\begin{equation*}
H_{n,\lambda}(A)=\mathbb{P}_n(A)-\lambda \mu(A),\ \mbox{for\ } A\in {\mathcal A}.\label{H1}
\end{equation*}

\begin{proposition}\label{polonik2} Let ${\mathcal A}$ be the class of compact sets $A$ with $ \textnormal{reach}(A)\geq r$ included in a given ball $B(0,R)$. Given a sample $X_1,\ldots,X_n$ from an absolutely continuous distribution $P$ with density $f$ in ${\mathbb R}^d$, let $S_n(\lambda)$ denote the empirical level set estimator
defined by $S_n(\lambda)=\mbox{\textnormal{argmin}}_{A\in{\mathcal A}}H_{n,\lambda}(A)$.
Assume that the level set $S(\lambda)=\{f\geq \lambda\}$
belongs to ${\mathcal A}$.
Then,

(a) The estimator $S_n(\lambda)$ is fully consistent to $S(\lambda)$.

(b) $L(S_n(\lambda))\to L(S(\lambda))$, a.s., where $L(A)$ denotes the outer Minkowski content of $A$.

\end{proposition}

\begin{proof}  (a) For simplicity, denote $H_{n,\lambda}=H_n$, $H_{\lambda}=H$, $S_n(\lambda)=S_n$ and $S(\lambda)=S$. First, note that, by definition, $H$ is a $d_\mu$-continuous functional and, from Theorem \ref{dHdmu}, it is also continuous with respect to $d_H$.
We have, with probability one,
{\small{
\begin{equation}
|\sup_{A\in{\mathcal A}}H_n(A)-\sup_{A\in{\mathcal A}}H(A)|\leq \sup_{A\in{\mathcal A}}|H_n(A)-H(A)|\leq \sup_{A\in{\mathcal A}}|{\mathbb P}_n(A)-P(A)|\rightarrow 0\label{GCC}
\end{equation}}}
since, from Theorem \ref{reachcompact}, ${\mathcal A}$ is a $P$-uniformity class.

From Theorem \ref{partialSn} and Proposition \ref{rconvexroll} to prove the full consistency we only need to establish
the $d_H$-convergence
\begin{equation}
d_H\left(S_n,S\right)\to 0,\ \mbox{a.s.}\label{dHSn}
\end{equation}
Now, let us take a value $\omega\in\Omega$ ($\Omega$ is the common probability space in which the random variables $X_i$ are defined) for which $d_H\left(S_n(\omega),S\right)\to 0$ does not hold but (\ref{GCC}) holds.
Then, as ${\mathcal A}$ is compact, we should have $d_H\left(S_n(\omega),T\right)\to 0$ for some subsequence $\{S_n(\omega)\}$ of $S_n=S_n(\omega)$ and for some $T\neq S$, $T\in {\mathcal A}$. Then, (\ref{GCC}) and the continuity of $H$ implies that
$H(T)=H(S)$. However, the convergence to $T\neq S$ is not possible, since $S\in{\mathcal A}$ is the (unique) maximum of $H$ in ${\mathcal A}$. Thus  we should have $d_H\left(S_n(\omega),S\right)\to 0$ for all $\omega$ such that (\ref{GCC}) holds. We conclude that (\ref{dHSn})  must hold with probability one.

(b) Follows again directly as a consequence of (a) together with Remarks 4.14 and 5.8 and Theorem 5.9 in Federer \cite{federer:59}.
\end{proof}

\

\noindent Our main point in this subsection is to show that our $P$-uniformity results fit in the framework developed by Polonik \cite{polonik:95}. This author obtains results of consistency (uniform in $\lambda$) of level set estimators, of type
 $\sup_\lambda d_\mu(S_n(\lambda),S(\lambda))\to 0$, a.s.
  One of his key assumptions is that the involved classes ${\mathcal C}$ are Glivenko-Cantelli classes.
To be more specific, the main result in Polonik's paper (which is Theorem 3.2) can be applied to our new $P$-uniformity class of sets with reach $\geq r$.

Let us first briefly recall the basic assumptions for the result in Polonik \cite{polonik:95}.

\begin{itemize}
\item[(A1)] For all $\lambda\geq 0$ the excess mass functional $H(C)=\int_C fd\mu-\lambda \mu(C)$ and its empirical counterpart $H_n(C)={\mathbb P}_n(C)-\lambda\mu(C)$ attain their maximum values on the class ${\mathcal C}$ at some sets denoted by $\Gamma(\lambda)$ and $\Gamma_n(\lambda)$, respectively.
\item[(A2)] The underlying density $f$ is bounded in ${\mathbb R}^d$.
\item[(A3)] ${\mathcal C}$ is a GC-class of closed sets with $\emptyset\in {\mathcal C}$.
\end{itemize}

Note that the maximizer over ${\mathcal C}$,  $\Gamma(\lambda)$, of $H(C)$ will only coincide with the level set
$S(\lambda)=\{f\geq \lambda\}$ whenever $S(\lambda)\in{\mathcal C}$. For this reason $\Gamma(\lambda)$ is called the \sl generalized $\lambda$-cluster\rm.

Now, let us denote by $g^*$ the ``measurable cover'' of any function $g$ (which of course coincides with $g$ when it is measurable). The main result in Polonik's paper is as follows:

\begin{theorem}\label{polonik} (Polonik (1995, Th. 3.2)) Let $\Lambda\subset[0,\infty)$. Suppose that, in addition to (A1)-(A3), the following two conditions hold:
\begin{itemize}
\item[(i)] For a distribution $Q$ in ${\mathbb R}^d$ with strictly positive density, the space $({\mathcal C},d_Q)$ is compact.
\item[(ii)] For every $\lambda\in\Lambda$ the generalized $\lambda$-cluster $\Gamma(\lambda)$ is unique up to $P$-null sets.
\end{itemize}
Then, with probability one,
\begin{equation}
\sup_{\lambda\in \Lambda}d_P(\Gamma(\lambda),\Gamma_n(\lambda))^*\rightarrow 0.\label{conclusion}
\end{equation}
\end{theorem}

Then, our point is that this result applies directly to the case in which ${\mathcal C}$ is the class ${\mathcal A}$ considered in Theorem \ref{reachcompact}, plus the empty set.
This is made explicit in the following statement.

\begin{theorem}\label{polonik-bis}
Under assumptions(A1), (A2) and (ii) of Theorem \ref{polonik}, the conclusion (\ref{conclusion}) is valid if we take ${\mathcal C}={\mathcal A}\cup\{\emptyset\}$, where ${\mathcal A}$ is the class of compact subsets of a given ball $B(0,R)$ with  \textnormal{reach} $\geq r$.
\end{theorem}

\begin{proof} We only need to prove that conditions (i) and (A3) hold in this case. The validity of (i) follows easily since ${\mathcal A}$ is $d_H$-compact from Theorems \ref{partialSn}  and \ref{reachcompact}, and therefore, from Theorem \ref{dHdmu}, it is also $d_Q$-compact for any absolutely continuous distribution $Q$. Also ${\mathcal C}={\mathcal A}\cup \{\emptyset\}$ is $d_Q$-compact.

As for the condition (A3) it also holds as a direct consequence of Theorem \ref{reachcompact}.
\end{proof}

Note that none of the typical examples of GC-classes (convex sets, ellipsoids, balls,...) allows us to consider multimodal densities. In this sense, Theorem \ref{polonik-bis} can be seen as an extension of uniformity results in the excess mass approach
beyond the realm of convex sets. The class of sets with reach bounded from below is much larger and includes non-connected members which are now candidates to be considered as possible level sets for multimodal densities.

\

\section*{Appendix: some numerical comparisons}\label{appendix}

In Theorem \ref{SSS} we have proved that the boundary length $L(S)$ of an $r$-convex set $S\subset{\mathbb R}^2$ can be consistently estimated in a plug-in way by $L(S_n)$, where $S_n$ is the $r$-convex hull of the sample.

Another estimator of $L(S)$ (in fact of $L_0(S)$) has been recently proposed by Jim\'{e}nez and Yukich \cite{jimenez:10}, based on the use of Delaunay triangulations. This estimator does not rely on any assumption of $r$-convexity but requires the use of sample data inside and outside the set $S$. The numerical results reported by Jim\'{e}nez and Yukich \cite{jimenez:10} show a remarkable performance of their \sl sewing based estimator\/ \rm (denoted by $L_n^s(S)$) which in fact outperforms that proposed in Cuevas et al. \cite{cuevas:07}.

Our plug-in estimator is not directly comparable with $L_n^s(S)$ since the required sampling models are different in both cases. Moreover, $L(S_n)$ incorporates the shape assumption of $r$-convexity on the target set $S$. Nevertheless, it is still interesting analyzing to what extent the use of the $r$-convexity assumption in $L(S_n)$ could improve the efficiency in the estimation.

\begin{figure}[!h]
\begin{center}
\setlength{\unitlength}{1mm}
\begin{picture}(120,50)
\put(0,-5){\scalebox{0.3}{\includegraphics{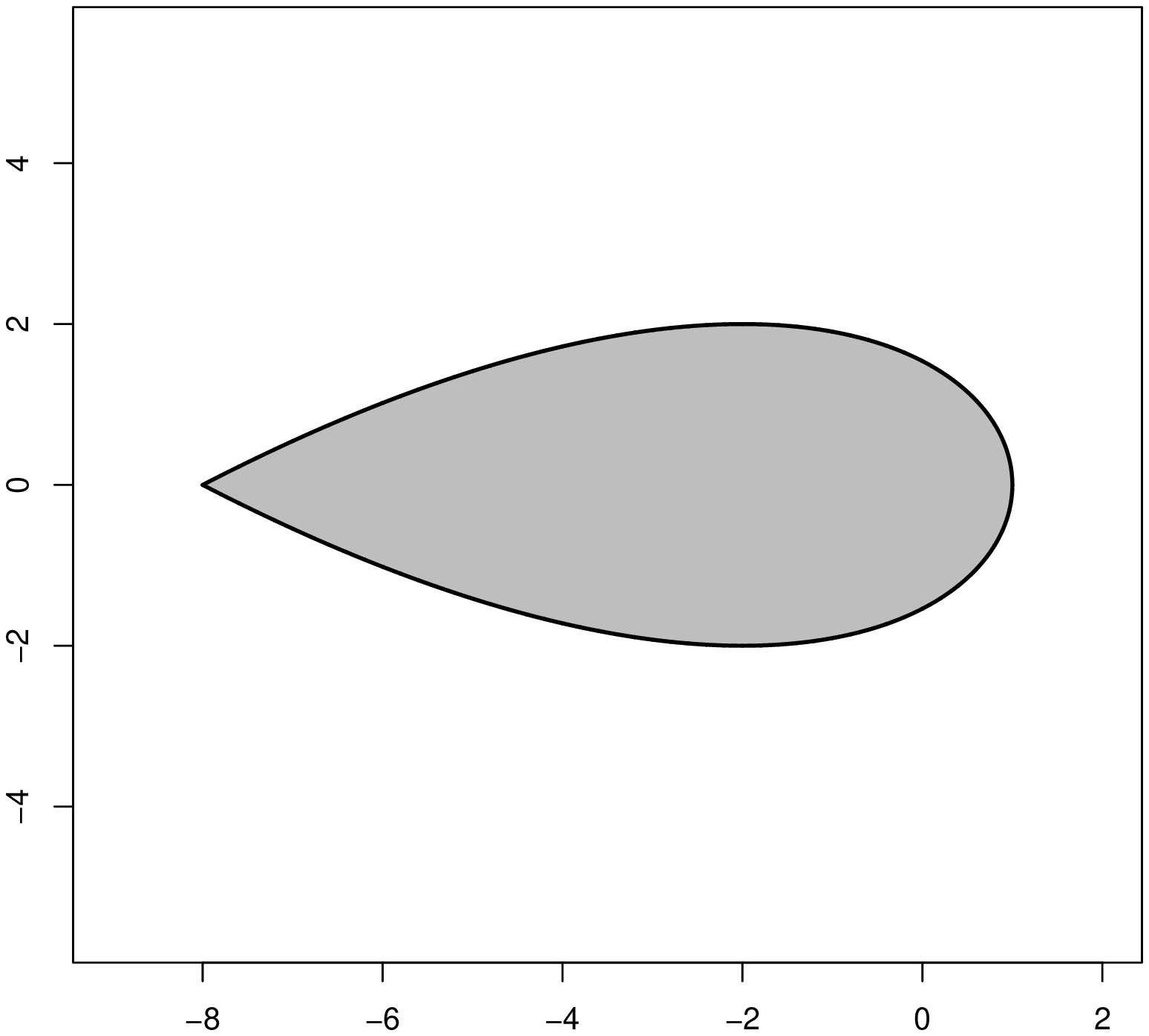}}}
\put(60,-5){\scalebox{0.3}{\includegraphics{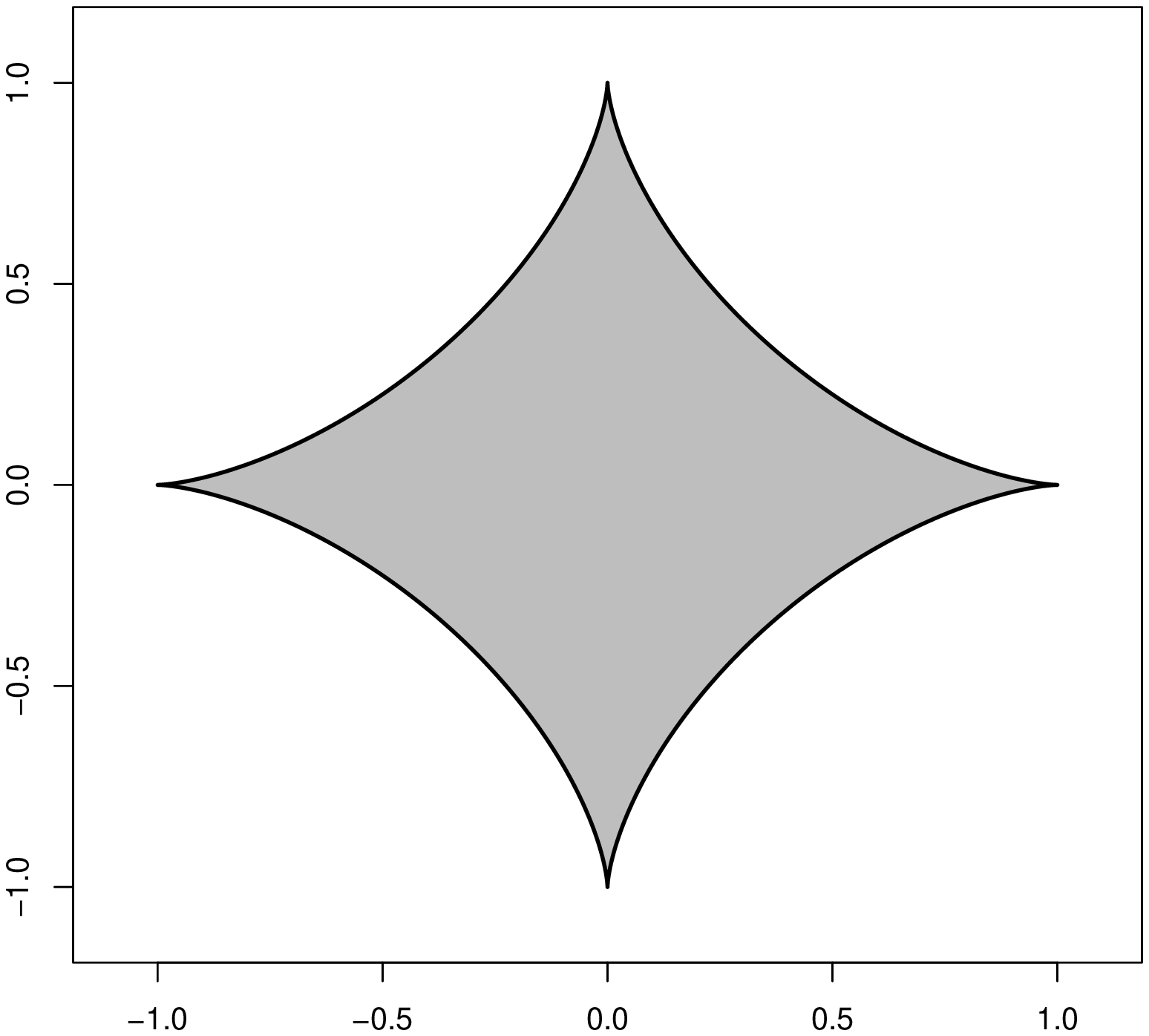}}}
\end{picture}
\caption{The supports in the simulation study are the domains inside a Catalan's trisectrix (left) and an astroid (right).}
\label{fig:curve}
\end{center}
\end{figure}

We have checked this, through a small simulation study. We have considered two $r$-convex sets $S$, defined as the domains inside two well-known closed curves: the \sl Catalan's trisectrix \/ \rm (as in Jim\'{e}nez and Yukich \cite{jimenez:10}) and the \sl astroid\rm. See Figure \ref{fig:curve}. The ``true'' maximal value of $r$ in the first set is $r=\infty$ (since it is convex). For the second set $r$ is close to 1. In practice, these values are not known so one must assume them as a model hypothesis keeping in mind that small values of $r$ correspond to more conservative (safer) choices. Recall that the class of $r$-convex sets is increasing as $r$ decreases.

The considered sample sizes are $n=5000,10000$. In the case of the estimator in Jim\'{e}nez and Yukich \cite{jimenez:10}, the sample observations are uniformly generated on a square containing the domain $S$ and it is assumed that we know (without error) whether a sample point belongs to $S$ or not. In the case of our estimator $L(S_n)$ all the observations are drawn from a uniform distribution on $S$.
The results are summarized in Table \ref{tab:curve}. The reported values in the columns ``mean'' and ``std'' correspond to averages over 1000 runs. The true values of $L(S)$ are 20.7846 (for the Catalan's trisectrix) and 6 (for the astroid). The outputs show a better behavior of $L(S_n)$, especially in terms of variability.

\begin{table}[!h]
\begin{center}
\caption{Mean and standard deviation of the sewing-based estimator $L_n^s(S)$ and the plug-in estimator $L(S_n)$, where $S_n$ is the $r$-convex hull of the sample. The estimator $S_n$ is computed for different values of $r$. The reported values correspond to averages over 1000 runs of samples size $n=5000$ and $n=10000$.}
\begin{tabular}{l|llcccc}
\multicolumn{3}{c}{}&\multicolumn{2}{c}{$n$=5000}&\multicolumn{2}{c}{$n$=10000}\\
\cline{4-7}
\multicolumn{3}{c}{}&Mean&Std&Mean&Std\\
  \hline
$S$ Catalan's trisectrix&$L_n^s(S)$& & 20.6684 & 0.3691 & 20.6966 & 0.3051\\
     $L(S)=20.7846$&$L(S_n)$&$r=0.5$ & 20.8722 & 0.0697 & 20.8289 & 0.0479\\
                            & &$r=2$ & 20.6383 & 0.0637 & 20.6821 & 0.0454\\
                             &&$r=5$ & 20.6131 & 0.0633 & 20.6661 & 0.0453\\
   \hline
       $S$ astroid&$L_n^s(S)$& & 5.5554 & 0.1324 & 5.6393 & 0.1061\\
  $L(S)=6$&$L(S_n)$&$r=0.25$ & 5.7024 & 0.0695 & 5.7709 & 0.0539\\
                  & &$r=0.5$ & 5.6848 & 0.0705 & 5.7584 & 0.0542\\
                     &&$r=1$ & 5.6730 & 0.0698 & 5.7500 & 0.0537\\
     \hline
\end{tabular}
\label{tab:curve}
\end{center}
\end{table}

\section*{Acknowledgements}
This work has been partially supported by Spanish Grants MTM2010-17366 and CCG10-UAM/ESP-5494 (A. Cuevas), MTM2010-17366 and Argentinian grant PIC-2008-0921 (R. Fraiman)
and by Spanish Grant MTM2008-03010 and the IAP research network grant no. P6/03 from the Belgian government (B. Pateiro-L\'{o}pez). The authors are grateful to Julio Gonz\'alez-D\'{\i}az for useful conversations.

\end{document}